\documentclass[11pt]{article}
\usepackage{amsmath,amssymb,amsthm}

\newcommand{\prf}{\noindent{\bf Proof.}\ }
\newcommand{\Tbar}{\ensuremath{\overline{\mathcal{T}}}}

\newcommand{\hess}{\operatorname{Hess}}
\newcommand{\grad}{\operatorname{grad}}

\newcommand{\caF}{\mathcal{F}}

\newcommand{\caH}{\mathcal{H}}

\newcommand{\caT}{\mathcal{T}}
\newcommand{\DR}{R_0\cup_{\mathcal N}\overline{R_0}}

\newtheorem{definition}{Definition}[section]
\newtheorem{lemma}[definition]{Lemma}
\newtheorem{theorem}[definition]{Theorem}

\newtheorem*{theorem'}{Theorem}
\newtheorem*{corollary'}{Corollary}
\newtheorem*{lemma'}{Lemma}

\begin{document}

\title{Extension of the Weil-Petersson connection\footnote{2000 Mathematics Subject Classification Primary: 32G15; Secondary: 20H10, 30F60.}}         
\author{Scott A. Wolpert}        
\date{September 16, 2007}          
\maketitle

\begin{abstract}
Convexity properties of Weil-Petersson geodesics on the Teichm\"{u}ller space of punctured Riemann surfaces are investigated.  A normal form is presented for the Weil-Petersson Levi-Civita connection for pinched hyperbolic metrics.  The normal form is used to establish approximation of geodesics in boundary spaces.  Considerations are combined to establish convexity along Weil-Petersson geodesics of the functions the distance between horocycles for a hyperbolic metric.
\end{abstract}

\section{Introduction}
Let $\caT$ be the Teichm\"{u}ller space of marked genus $g$, $n$ punctured Riemann surfaces with hyperbolic metrics.  Associated to hyperbolic metrics on Riemann surfaces are quantities on Teichm\"{u}ller space: the Weil-Petersson (WP) metric, the geodesic-length functions and for $n>0$ the functions with values the distances between unit-length horocycles. A fundamental property of geodesic-length functions is their convexity along WP geodesics. We now generalize the property and show that the distance between horocycles is also convex along WP geodesics.  In particular we show for surfaces with punctures that the distance between horocycles and WP geodesics are suitably approximated by geodesic-lengths and WP geodesics for surfaces without punctures.  The convexity of the distance between horocycles provides new information on the behavior of WP geodesics.

There are applications of the convexity of geodesic-length functions.  The convexity and properness for a set of {\em filling} geodesic-lengths provides for a solution of the Nielsen realization problem and provides that with the WP metric $\mathcal T$ is a convex space \cite{Wlnielsen}.  The convexity and refined forms of Masur's expansion are the essential ingredients in showing that the augmented Teichm\"{u}ller space is $CAT(0)$ with all strata convex 
\cite{DW2,Wlcomp,Yam}.  The convexity also provides for the convexity of Bers regions, an important consideration in establishing Brock's quasi isometry between $\caT$ and the pants graph \cite{Brkwp}.

The decorated Teichm\"{u}ller space $\mathcal D$ for surfaces with cusps is a principal $\mathbb R_+^n$-bundle over $\caT$ \cite{Pendec, Pencell}.  Penner continues his investigation of the mapping class group invariant geometry of $\mathcal D$ and $\caT$.  Mapping class group invariant cell decompositions are obtained from decompositions of hyperbolic surfaces involving disjoint collections of simple geodesics between cusps.  The cells are indexed by ribbon graphs (fatgraphs) and the decompositions are considerations for several celebrated results.  Maximal collections of disjoint simple geodesics between cusps also give rise to global coordinates.  Considerations begin with Penner's {\em lambda length} parameter, the root exponential distance between horocycles $\lambda=(2e^{\delta})^{1/2}$ \cite{Pendec, Pencell}. The description of $\caT$, the WP symplectic form, and the action of mapping classes is particularly straightforward.  Most recently McShane and Penner investigate the degeneration of hyperbolic metrics and determine which geodesic-lengths can be small in a given cell \cite{McP}.  Also recently Mondello considers the Teichm\"{u}ller space $\caT_b$ of hyperbolic surfaces with geodesic boundaries and investigates the distance between boundaries as a generalization of the distance between horocycles.  He finds that the WP Poisson geometry of $\caT_b$ limits to a Poisson geometry of 
$\caT$ \cite{Mon}. 

Masur first considered the expansion of the WP metric for surfaces with short geodesics (pinched surfaces) 
\cite{Msext}.  Refined expansions for the metric are presented in \cite{Wlcomp,Wlbhv}.  An expansion for the WP Levi-Civita connection is required for a detailed understanding of WP geodesics for surfaces with short geodesics.  In Section \ref{conWP} we develop and present such an expansion in Theorem \ref{WPconn}.  The expansion is a normal form in terms of the short geodesic-length gradients $\{\grad \ell_{\alpha}^{1/2}\}_{\alpha\in\sigma}$ and a relative length basis $\{\grad\ell_{\beta}\}$ for the remaining gradients of geodesic-lengths.  To properly analyze surfaces with short geodesic-lengths the augmented Teichm\"{u}ller space $\Tbar$, a partial compactification, is introduced.  The geometry in a neighborhood of a boundary space $\caT(\sigma)=\{\ell_{\alpha}=0\mid\alpha\in\sigma\}\subset\Tbar$ of marked degenerate hyperbolic metrics is considered.  In Section \ref{WPgeoapp} we develop and present an appropriate $C^1$-approximation of geodesics in $\caT(\sigma)$ by geodesics in $\caT$ with $\ell_{\alpha},\alpha\in\sigma,$ small.  The approximation is presented for a geodesic $\gamma(t)$ with the initial values  
$\ell_{\alpha}(\gamma(0)),$ $\langle\grad\ell_{\alpha}^{1/2},\dot\gamma(0)\rangle,$ $\langle J\grad\ell_{\alpha}^{1/2},\dot\gamma(0)\rangle,$ $\alpha\in\sigma$ all small and $J$ the almost complex structure.  From an expansion for the Hessian of $\ell_{\alpha}$ and Theorem \ref{approx} the geodesic-lengths $\ell_{\alpha}$ on $\gamma(t)$ are $C^2$-almost constant.  There are geodesics nearly parallel to the boundary spaces of $\Tbar$.    

Distance between horocycles is also described in terms of lengths of closed geodesics for a degenerating family of hyperbolic metrics and degenerating collar widths $2w(\alpha)$.   For a hyperbolic geodesic $\omega_0$ between cusps $p,q$ for a surface $R_0$ we introduce the double across the cusps $R_0\cup_{p,q}\overline{R_0}$ and open the cusps $p,q$ to obtain an approximating surface $R$ with an approximating closed geodesic $\omega$.  The closed geodesic $\omega$ consists of segments in the thick regions of $R$ and segments crossing the collars $\mathcal C(p)$ and $\mathcal C(q)$. The length $\ell_{\omega}-2w(p)-2w(q)$ of $\omega$ with the collar widths subtracted approximates twice the distance 
$2\tilde\ell_{\omega_0}$ between horocycles.  We use the description of the distance $\tilde\ell_{\omega_0}$ in terms of closed geodesics and collar widths to establish the main result in Theorem \ref{horowpcon}.  Our approach suggests that the convexity of $\ell_{\omega}$ approximates the convexity of $\ell_{\omega_0}$. 

The discussion is arranged as follows.  Certain aspects of hyperbolic geometry are reviewed in Section \ref{basics}.  Formulas for the WP gradient and Hessian of geodesic-length are presented in Section \ref{gradhess}.  A normal form for the WP metric and covariant derivative are developed in Section \ref{conWP}.  The boundary spaces and approximation of geodesics of $\Tbar$ is described in Section \ref{WPgeoapp} and the main result is developed in Section \ref{horocon}.     

The author thanks Yair Minsky for raising the question of convexity of the distance between horocycles.         

\section{Collars, cusp regions and the mean value estimate}
\label{basics}
A Riemann surface with hyperbolic metric can be considered as the union of a $thick$ region where the injectivity radius is bounded below by a positive constant and a complementary $thin$ region.  The totality of all $thick$ regions of Riemann surfaces of a given topological type forms a compact set of metric spaces in the Gromov-Hausdorff topology.  A $thin$ region is a disjoint union of collar and cusp regions.  We describe basic properties of collars and cusp regions including bounds for the injectivity radius and separation of simple geodesics. 

We follow Buser's presentation \cite[Chap. 4]{Busbook}.   For a geodesic $\alpha$ of length $\ell_{\alpha}$ on a Riemann surface the collar about the geodesic is $\mathcal C(\alpha)=\{d(p,\alpha)\le w(\alpha)\}$ for the width $w(\alpha)$, $\sinh w(\alpha)\sinh \ell_{\alpha}/2=1$.  The width is given as $w(\alpha)=\log 4/\ell_{\alpha}+O(\ell_{\alpha}^2)$ for $\ell_{\alpha}$ small.  For $\mathbb H$ the upper half plane with hyperbolic distance $d(\ ,\ )$ a collar is covered by the region  $\{d(z,i\mathbb R^+)\le w(\alpha)\}\subset \mathbb H$ with deck transformations generated by $z\rightarrow e^{\ell_{\alpha}}z$.  The quotient $\{d(z,i\mathbb R^+)\le w(\alpha)\}\slash \langle z\rightarrow e^{\ell_{\alpha}}z\bigr >$ embeds into the Riemann surface. For $z$ in $\mathbb H$ the prescribed region is approximately $\{\ell_{\alpha}/2\le \arg z\le \pi-\ell_{\alpha}/2\}$.   A cusp region $\mathcal C_{\infty}$ is covered by the region $\{\Im z\ge 1/2\}\subset\mathbb H$ with deck transformations generated by $z\rightarrow z+1$.  The quotient $\{\Im z\ge 1/2\}\slash \langle z\rightarrow z+1\rangle $ embeds into the Riemann surface.  The boundary of a collar $\mathcal C(\alpha)$ for $\ell_{\alpha}$ bounded and boundary of a cusp region $\mathcal C_{\infty}$ have length approximately $2$.

\begin{theorem}
\label{collars}
For a Riemann surface of genus $g$ with $n$ punctures given pairwise disjoint simple closed geodesics $\alpha_1,\dots,\alpha_m$ there exist simple closed geodesics $\alpha_{m+1},\dots,\alpha_{3g-3+n}$ such that $\alpha_1,\dots,\alpha_{3g-3+n}$ are pairwise disjoint.  The collars $\mathcal C(\alpha_j)$ about $\alpha_j$, $1\le j\le 3g-3+n$, and the cusp regions are mutually pairwise disjoint.
\end{theorem}

On $thin$ the injectivity radius is bounded below in terms of the distance into a collar or cusp region.  For a point 
$p$ of a collar or cusp region of a Riemann surface write $inj(p)$ for the injectivity radius and $\delta(p)$ for the distance to the boundary of the collar  or cusp region.   The injectivity radius is bounded as follows, \cite[II, Lemma 2.1]{Wlspeclim}.

\begin{lemma}
\label{enhcollar}
The product $inj(p)\,e^{\delta(p)}$ of injectivity radius and exponential distance to the boundary is bounded below by a positive constant.  
\end{lemma}
\noindent The standard consideration for simple closed geodesics and cusp regions generalizes as follows, \cite[Lemma 2.3]{Wlbhv}.

\begin{lemma}
\label{separ}
A simple closed geodesic is disjoint from the cusp regions.   A simple closed geodesic is either disjoint from a collar or crosses the collar or is the core geodesic. 
\end{lemma}

We use the following mean value estimate for functions on the upper half plane with hyperbolic area element $dA$, \cite{Fay}.

\begin{lemma}
\label{enhmv}
Harmonic Beltrami differentials and the exponential-distance function $e^{-2d(p,\tilde\alpha)}$ for a geodesic $\tilde\alpha$ satisfy a mean value estimate on $\mathbb H$ with constant determined by the radius of the ball.  A non negative function $f$ satisfying a mean value estimate 
on $\mathbb H$ and a subset of a discrete group $\mathcal G\subset\Gamma$ satisfy
\[
\sum_{A\in\mathcal G} f(A(p))\le c\,inj(p)^{-1}\int_{\cup_{A\in \mathcal G}A(\mathbf B(p,\epsilon))}f\,dA
\]
with constant determined by the mean value constant.
\end{lemma}

\section{The gradient and Hessian of geodesic-length}
\label{gradhess}

We review properties of the gradient and Hessian of geodesic-length as presented in \cite{Wlbhv}.  The reader should consult the reference for proofs and further applications. 

Points of the Teichm\"{u}ller space $\caT$ are equivalence classes $\{(R,ds^2,f)\}$ of marked complete hyperbolic structures with reference homeomorphisms $f:F\rightarrow R$ from a base surface $F$.  Basic invariants of a hyperbolic metric are the lengths of the unique closed geodesic representatives of the non peripheral free homotopy classes.  For the non peripheral free homotopy class $\alpha$ on $F$ the length of the geodesic representative for $f(\alpha)$ is the value of the geodesic-length $\ell_{\alpha}$ at the marked structure.  For $R$ with uniformization representation $f_*:\pi_1(F)\rightarrow\Gamma \subset PSL(2;\mathbb R)$ and $\alpha$ corresponding to the conjugacy class of an element $A$ then $\cosh \ell_{\alpha}/2= tr A/2$.  Collections of geodesic-lengths provide local coordinates for $\caT$, \cite{Busbook,ImTan,WlFN}. 

From Kodaira-Spencer deformation theory the infinitesimal deformations of a surface $R$ are represented by the Beltrami differentials $\caH (R)$ harmonic with respect to the hyperbolic metric, \cite{Ahsome}. Also the cotangent space of $\caT$ at $R$ is $Q(R)$ the space of holomorphic quadratic differentials with at most simple poles at the punctures of $R$.  The holomorphic tangent-cotangent pairing is
\[
(\mu,\varphi)=\int_R \mu \varphi
\]
for $\mu\in\caH(R)$ and $\varphi\in Q(R)$.  Elements of $\caH (R)$ are symmetric tensors given as $\overline\varphi(ds^2)^{-1}$ for $\varphi\in Q(R)$ and $ds^2$ the hyperbolic metric.  The Weil-Petersson (WP) Hermitian metric and cometric are given as
\[
\langle\mu,\nu\rangle_{Herm}=\int_R\mu\overline\nu dA\quad\mbox{and}
\quad\langle\varphi,\psi\rangle_{Herm}=\int_R\varphi\overline\psi (ds^2)^{-1}
\]
for $\mu,\nu\in\caH (R)$ and $\varphi, \psi\in Q(R)$ and $dA$ the hyperbolic area element.  The WP Riemannian metric is $\langle\ ,\ \rangle =\Re\langle\ ,\ \rangle_{Herm}$.  The WP  metric is K\"{a}hler, non complete, with non pinched negative sectional curvature and determines a $CAT(0)$ geometry, see \cite{Ahsome, Zh2, Ngbook, Trmbook, Wlcomp} for references and background.

For a closed geodesic $\alpha$ conjugate the group $\Gamma$ for the geodesic to correspond to the deck transformation 
$A:z\rightarrow e^{\ell_{\alpha}}z$ with imaginary axis $\tilde\alpha$ and for $p$ in $\mathbb H$ with coordinate $z$ consider the pair of coset sums
\begin{equation}
\label{series}
\mathbb P_{\alpha}(p)=\sum_{B\in\langle A\rangle\backslash\Gamma}e^{-2d(B(p),\tilde\alpha)} \quad\mbox{and}\quad \Theta_{\alpha}=\sum_{B\in\langle A\rangle\backslash\Gamma} B^*\bigl(\frac{dz}{z}\bigr)^2
\end{equation}
for $\langle A\rangle$ the cyclic group generated by $A$.  For $\mathcal F$ a $\Gamma$ fundamental domain an unfolding
\[
\int_{\mathcal F}\sum_{B\in\langle A\rangle\backslash\Gamma}f\circ B\,dA\, = \sum_{B\in\langle A\rangle\backslash\Gamma}
\int_{B(\mathcal F)}f\,dA\,=\int_{\mathbb H/\langle A \rangle}f\,dA
\]
provides that $\overline{\Theta_{\alpha}}(ds^2)^{-1}$ and $\mathbb P_{\alpha}$ are bounded in $L^1$ norm by a multiple
of $\ell_{\alpha}$ and by Lemma \ref{enhmv} that the coset sums are convergent.  
The WP gradient of geodesic-length is $\grad \ell_{\alpha}=\frac{2}{\pi}\overline{\Theta_{\alpha}}(ds^2)^{-1}$ with the infinitesimal Fenchel-Nielsen right-twist deformation given as $t_{\alpha}=\frac12J\grad\ell_{\alpha}$ for $J$ the almost complex structure of $\caT$, \cite{Gardtheta, WlFN}.  We combined the formula of Riera for pairing geodesic-length gradients with an analysis of the orbit of a discrete group to develop an expansion for the pairing, \cite{Wlbhv}.
\begin{theorem}
\label{gradpair}
The WP pairing of gradients for disjoint geodesics $\alpha,\beta$ satisfies
\[
\langle\grad\ell_{\alpha},\grad\ell_{\beta}\rangle =\frac{2}{\pi}\ell_{\alpha}\delta_{\alpha\beta}\,+\,O(\ell_{\alpha}^2\ell_{\beta}^2)
\]
where for $c_0$ positive the constant for the positive remainder term is uniform for 
$\ell_{\alpha},\ell_{\beta}\le c_0$.
\end{theorem}  

The Fenchel-Nielsen twist deformation and geodesic-length are related by duality in the WP K\"{a}hler form and satisfy
\begin{equation}
\label{cos}
\langle\grad \ell_{\alpha},J\grad\ell_{\beta}\rangle\,=\,-2\sum_{p\in\alpha\cap\beta}\cos \theta_p
\end{equation}
for the sum over transverse intersections.  The sum vanishes if the geodesics are disjoint or coincide.  For the coordinate $z$ of $\mathbb H$ given in polar form $re^{i\theta}$ the hyperbolic distance to the imaginary axis satisfies
$e^{-2d(z,i\mathbb R^+)}\le \sin^2\theta \le 4 e^{-2d(z,i\mathbb R^+)}$.   The inequalities 
$| B^*(\frac{dz}{z})^2(ds^2)^{-1}|\le 4e^{-2d(B,\tilde\alpha)}$ and  $|\Theta_{\alpha}(ds^2)^{-1}|\le 4\mathbb P_{\alpha}$ are consequences.

We review the description of the WP Hessian of geodesic-length, \cite{Wlbhv}.  It is natural to introduce a model for the hyperbolic plane presenting the transformation $z\rightarrow\lambda z$ as a translation.  The horizontal strip $\mathbb S=\{\zeta\mid 0<\Im\zeta<\pi\}$ is a suitable model.  The upper half plane and horizontal strip are related by the map $\zeta=\log z$.  We consider for a Beltrami differential $\mu\in\mathcal H(\Gamma)$ and $t$ small the family of qc homeomorphisms $g^t$ of $\mathbb S$ satisfying $g^t_{\overline{\zeta}}=t\mu g^t_{\zeta}$ and conjugating the translation $\zeta\rightarrow\zeta + \ell$ to a translation $\zeta\rightarrow\zeta+\ell^t$, \cite{AB}.  We have from the translation equivariance the equation
\[
g^t(\zeta+\ell)=g^t(\zeta)+\ell^t.
\]
The derivatives of geodesic-length are given for the fundamental domain $\mathcal F=\{\zeta\in\mathbb S\mid 0\le\Re\zeta<\ell\}$ and deformation vector field $\dot{g}=\frac{d}{dt}g^t\mid_{t=0}$ as follows.  
\begin{lemma}
\label{ldot}
For $\mu\in\mathcal H(\Gamma)$ with deformation vector field $\dot{g}$ then
\[
\dot{\ell}=\frac{2}{\pi}\Re\int_{\mathcal F}\mu\,dE\quad\mbox{and}\quad\ddot{\ell}=\frac{4}{\pi}\Re\int_{\mathcal F}\mu\dot{g}_{\zeta}\,dE
\]
for $dE$ the Euclidean area element.
\end{lemma}

The series $\mathbb P_{\alpha}$ serves as a weight function to bound the convexity of geodesic-length functions, \cite{Wlbhv}.
\begin{theorem}
\label{hessbd}
The Hessian of geodesic-length is bounded in terms of the weight function $\mathbb P_{\alpha}$ for the WP pairing
\[
\langle\mu,\mu\mathbb P_{\alpha}\rangle \le 3\pi\hess\ell_{\alpha}[\mu,\mu] \le 48\langle\mu,\mu\mathbb P_{\alpha}\rangle
\]
for $\mu\in\mathcal H(\Gamma)$.
\end{theorem}

The root geodesic-length function $\ell_{\alpha}^{1/2}$ is also convex with positive definite Hessian.  We review the basic expansion for the WP connection $D$ and the root geodesic-length gradient $\lambda_{\alpha}=\grad\ell_{\alpha}^{1/2}$, \cite{Wlbhv}.  The Riemannian Hessian is directly related to covariant differentiation $\hess h(U,V)=\langle D_U\grad h,V\rangle$ for vector fields $U,V$, \cite{BO'N}.

\begin{theorem}
\label{conn}
The WP connection $D$ for the root geodesic-length gradient $\lambda_{\alpha}$ satisfies
\[
D_U\lambda_{\alpha}=3\ell_{\alpha}^{-1/2}\langle J\lambda_{\alpha},U\rangle J\lambda_{\alpha}+O(\ell_{\alpha}^{3/2}\|U\|)
\]
in terms of the WP norm $\|\ \|$ where for $c_0$ positive the remainder term constant is uniform for 
$\ell_{\alpha}\le c_0.$
\end{theorem}

\section{Continuity of the WP pairing and connection}
\label{conWP}

The WP completion is the augmented Teichm\"{u}ller space $\Tbar$.  The space $\Tbar$ is described
in terms of Fenchel-Nielsen coordinates and in terms of the Chabauty topology for $PSL(2;\mathbb R)$ representations.  The strata of $\Tbar$ correspond to collections of {\em vanishing lengths} and describe Riemann surfaces with nodes in the sense of Bers, \cite{Bersdeg}.  We present a local frame for the tangent bundle $T\mathcal T$ in terms of gradients of geodesic-length functions for a neighborhood of a stratum point.   The behavior of the WP pairing for the frame is presented in Lemma \ref{paircont} and of the WP connection for the frame is presented in Theorem \ref{WPconn}.  

The points of the Teichm\"{u}ller space $\mathcal T$ are equivalence classes 
$\{(R,f)\}$ of surfaces with reference homeomorphisms $f:F\rightarrow R$.  
The {\em complex of curves} $C(F)$ is defined as follows.  The vertices of
$C(F)$ are the free homotopy classes of homotopically nontrivial, non peripheral,
simple closed curves on $F$.  An edge of the complex consists of a pair of
homotopy classes of disjoint simple closed curves.  A $k$-simplex consists of
$k+1$ homotopy classes of mutually disjoint simple closed curves.  A maximal
set of mutually disjoint simple closed curves, a {\em partition},
has $3g-3+n$ elements.  The mapping class group $Mod$ acts on the complex $C(F)$.

The Fenchel-Nielsen coordinates for $\mathcal T$ are given in terms of geodesic-lengths and
lengths of auxiliary geodesic segments, \cite{Abbook,Busbook,ImTan}.  A partition
 $\mathcal{P}=\{\alpha_1,\dots,\alpha_{3g-3+n}\}$ decomposes the
reference surface $F$ into $2g-2+n$ components, each homeomorphic to
a sphere with a combination of three discs or points removed.  A marked
Riemann surface $(R,f)$ is likewise decomposed into pants by the geodesics
representing the elements of $\mathcal P$.  Each component pants, relative to its hyperbolic
metric, has a combination of three geodesic boundaries and cusps.  For each
component pants the shortest geodesic segments connecting boundaries determine
 designated points  on each boundary.  For each geodesic $\alpha$ in the pants
decomposition of $R$ a parameter $\tau_{\alpha}$ is defined as the displacement along
the geodesic between designated points, one  for each side of the geodesic.
For marked Riemann surfaces close to an initial reference marked Riemann surface, the
displacement $\tau_{\alpha}$ is the distance between the designated points; in
general the displacement is the analytic continuation (the lifting) of the
distance measurement.  For $\alpha$ in $\mathcal P$ define the {\em
Fenchel-Nielsen  angle} by $\vartheta_{\alpha}=2\pi\tau_{\alpha}/\ell_{\alpha}$.
The Fenchel-Nielsen coordinates for Teichm\"{u}ller space for the
decomposition $\mathcal P$ are
$(\ell_{\alpha_1},\vartheta_{\alpha_1},\dots,\ell_{\alpha_{3g-3+n}},\vartheta_{\alpha_{3g-3+n}})$.
The coordinates provide a real analytic equivalence of $\mathcal T$ to
$(\mathbb{R}_+\times \mathbb{R})^{3g-3+n}$, \cite{Abbook,Busbook,ImTan}.

A partial compactification of Teichm\"{u}ller space is introduced by extending the range
of the Fenchel-Nielsen parameters.  The added points correspond to unions of
hyperbolic surfaces with formal pairings of cusps.  The interpretation of {\em length
vanishing} is the key ingredient.  For an $\ell_{\alpha}$ equal to zero, the
angle $\vartheta_{\alpha}$ is not defined and in place of the geodesic for
$\alpha$ there appears a pair of cusps; the reference map $f$ is now a homeomorphism of
$F-\alpha$ to a union  of hyperbolic surfaces (curves parallel to $\alpha$
map to loops encircling the cusps).  The parameter space for a pair
$(\ell_{\alpha},\vartheta_{\alpha})$ will be the identification space $\mathbb{R}_{\ge
0}\times\mathbb{R}/\{(0,y)\sim(0,y')\}$.  More generally for the pants decomposition
$\mathcal P$ a frontier set $\mathcal{F}_{\mathcal P}$ is added to the
Teichm\"{u}ller space by extending the Fenchel-Nielsen parameter ranges: for
each $\alpha\in\mathcal{P}$, extend the range of $\ell_{\alpha}$ to include
the value $0$, with $\vartheta_{\alpha}$ not defined for $\ell_{\alpha}=0$.  The
points of $\mathcal{F}_{\mathcal P}$ parameterize unions of Riemann
surfaces with each $\ell_{\alpha}=0,\alpha\in\mathcal{P},$ specifying a pair
of cusps. The points of $\mathcal F_{\mathcal P}$ are Riemann surfaces with nodes in the sense of Bers, \cite{Bersdeg}.   For a simplex $\sigma\subset\mathcal{P}$, define the
$\sigma$-null stratum, a subset of $ \mathcal F_{\mathcal P}$, as $\mathcal{T}(\sigma)=\{R\mid \ell_{\alpha}(R)=0
\mbox{ iff }\alpha\in\sigma\}$. Null strata are given as products of lower dimensional
Teichm\"{u}ller spaces.   The frontier set $\mathcal{F}_{\mathcal P}$
is the union of the $\sigma$-null strata for the sub simplices of
$\mathcal{P}$.  Neighborhood bases for points of $\mathcal{F}_{\mathcal P}
\subset\mathcal{T}\cup\mathcal{F}_{\mathcal P}$  are specified by the
condition that for each simplex $\sigma\subset\mathcal P$  the projection
$((\ell_{\beta},\vartheta_{\beta}),\ell_{\alpha}):
\mathcal{T}\cup\mathcal{T}(\sigma)\rightarrow\prod_{\beta\notin\sigma}(\mathbb{R}_+\times\mathbb{R})\times\prod_{\alpha\in\sigma}(\mathbb{R}_{\ge
0})$  is continuous.    For a simplex $\sigma$ contained in partitions $\mathcal P$ and $\mathcal P'$ the specified neighborhood systems for $\mathcal T \cup\mathcal{T}(\sigma)$ are equivalent.  The {\em
augmented Teichm\"{u}ller space} $\Tbar=\mathcal{T}\cup_{\sigma\in
C(F)}\mathcal{T}(\sigma)$ is the resulting stratified topological space,
\cite{Abdegn, Bersdeg}.  $\Tbar$ is not locally compact since points of the
frontier do not have relatively compact neighborhoods; the  neighborhood bases are
unrestricted in the $\vartheta_{\alpha}$ parameters for $\alpha$ a $\sigma$-null.
The action of $Mod$ on $\mathcal T$ extends to an action by homeomorphisms on
$\Tbar$ (the action on $\Tbar$ is not properly discontinuous) and the quotient
$\Tbar/Mod$ is topologically the compactified moduli space of stable curves, \cite[see Math. Rev.56 \#679]{Abdegn}.

We present an alternate description of the frontier points in terms of representations of groups and the Chabauty topology.  A Riemann surface with punctures and hyperbolic metric is uniformized by a cofinite subgroup $\Gamma\subset PSL(2;\mathbb R)$.  A puncture corresponds to the $\Gamma$-conjugacy class of a maximal parabolic subgroup.  In general a Riemann surface with punctures corresponds to the $PSL(2;\mathbb R)$ conjugacy class of a tuple $(\Gamma,\langle\Gamma_{01}\rangle ,\dots,\langle\Gamma_{0n}\rangle )$ where $\langle\Gamma_{0j}\rangle $ are the maximal parabolic classes and a labeling for punctures is a labeling for conjugacy classes.  A {\em Riemann surface with nodes} $R'$ is a finite collection of $PSL(2;\mathbb R)$ conjugacy classes of tuples 
$(\Gamma^\ast,\langle\Gamma_{01}^\ast\rangle ,\dots,\langle\Gamma_{0n^\ast}^\ast\rangle )$ with a formal pairing of certain maximal parabolic classes.  The conjugacy class of a tuple is called a {\em part} of $R'$.  The unpaired maximal parabolic classes are the punctures of $R'$ and the genus of $R'$ is defined by the relation $Total\ area=2\pi(2g-2+n)$.  A cofinite $PSL(2;\mathbb R)$ injective representation of the fundamental group of a surface is topologically allowable provided peripheral elements correspond to peripheral elements.  A point of the Teichm\"{u}ller space $\mathcal T$ is given by the $PSL(2;\mathbb R)$ conjugacy class of a topologically allowable injective cofinite representation of the fundamental group $\pi_1(F)\rightarrow\Gamma\subset PSL(2;\mathbb R)$.  For a simplex $\sigma$ a point of $\mathcal T(\sigma)$ is given by a collection $\{(\Gamma^\ast,\langle\Gamma_{01}^\ast\rangle ,\dots,\langle\Gamma_{0n^\ast}^\ast\rangle )\}$ of tuples with: a bijection between $\sigma$ and the paired maximal parabolic classes; a bijection between the components $\{F_j\}$ of $F-\sigma$ and the conjugacy classes of parts $(\Gamma^j,\langle\Gamma_{01}^j\rangle ,\dots,\langle\Gamma_{0n^j}^j\rangle )$ and the $PSL(2;\mathbb R)$ conjugacy classes of topologically allowable isomorphisms $\pi_1(F_j)\rightarrow\Gamma^j$, \cite{Abdegn, Bersdeg}. We are interested in geodesic-lengths for a sequence of points of $\mathcal T$ converging to a point of $\mathcal T(\sigma)$.  The convergence of hyperbolic metrics provides that for closed curves of $F$ disjoint from $\sigma$ geodesic-lengths converge, while closed curves with essential $\sigma$ intersections have geodesic-lengths tending to infinity, \cite{Bersdeg, Wlhyp}.  

We refer to the Chabauty topology to describe the convergence for the $PSL(2;\mathbb R)$ representations.  
Chabauty introduced a topology for the space of discrete subgroups of a locally compact group, \cite{Chb}.   A neighborhood of $\Gamma\subset PSL(2;\mathbb R)$ is specified by a neighborhood $U$ of the identity in $PSL(2;\mathbb R)$ and a compact subset $K\subset PSL(2;\mathbb R)$.  A discrete group $\Gamma'$ is in the neighborhood $\mathcal N(\Gamma,U,K)$ provided $\Gamma'\cap K\subseteq\Gamma U$ and $\Gamma\cap K\subseteq\Gamma'U$.  The sets $\mathcal N(\Gamma,U,K)$ provide a neighborhood basis for the topology. We consider a sequence of points of $\mathcal T$ converging to a point of $\mathcal T(\sigma)$ corresponding to 
$\{(\Gamma^\ast,\langle\Gamma_{01}^\ast\rangle ,\dots,\langle\Gamma_{0n^\ast}^\ast\rangle )\}$.  Important for the present considerations is the following property.  Given a sequence of points of $\mathcal T$ converging to a point of $\mathcal T(\sigma)$ and a component $F_j$ of $F-\sigma$ there exist $PSL(2;\mathbb R)$ conjugations such that restricted to $\pi_1(F_j)$ the corresponding representations converge element wise to $\pi_1(F_j)\rightarrow\Gamma^j$, \cite[Thrm. 2]{HrCh}.

We now consider the WP geometry of geodesic-length functions in a neighborhood of an augmentation point $p$ of 
$\mathcal T(\sigma)\subset\overline{\mathcal T}$.  The null stratum $\caT(\sigma)$ is metrically characterized as the union of all WP distance-realizing paths containing $p$ as an interior point.  For the following considerations 
we refer to the elements of $\sigma$ as the {\em short geodesics}. 
We are interested in collections of geodesic-length gradients providing a frame for the 
tangent bundle in a neighborhood of $p$.
\begin{definition}
\label{relbas}
A relative length basis for a point $p$ of $\mathcal T(\sigma)$ is a collection $\tau$ of vertices of $C(F)$   disjoint from the elements of $\sigma$ such that at $p$ the gradients $\{\grad \ell_{\beta}\}_{\beta\in\tau}$ provide the germ of a frame over $\mathbb R$ for the tangent space $T\mathcal T(\sigma)$.
\end{definition}
\noindent Examples of relative length bases are given as the union of a partition and a dual partition 
for $R-\sigma$, see \cite[Chap. 3 Secs. 3 \& 4]{Busbook}, \cite[{\em markings} in Sec. 2.5]{MaMiII} and \cite[Thrm. 3.4]{WlFN}.   The lengths of the elements of a relative length basis are necessarily bounded on a neighborhood in $\overline{\mathcal T}$.  We consider below for a geodesic $\gamma$ the root geodesic-length $\ell^{1/2}_{\gamma}$ and gradient $\lambda_{\gamma}=\grad \ell^{1/2}_{\gamma}$.  We introduce the convention that on $\caT(\sigma)$ the pairing $\langle\grad\ell_{\alpha},\grad\ell_{\beta}\rangle $ vanishes for geodesics 
$\alpha,\beta$ on distinct components (parts) of a Riemann surface with nodes.  
\begin{lemma}
\label{paircont}
The WP pairing of geodesic-length gradients $\langle\grad \ell_{\alpha},\grad \ell_{\beta}\rangle $ is continuous in a neighborhood of a point $p$ of $\mathcal T(\sigma)\subset\overline{\mathcal T}$ for $\alpha$ and $\beta$ disjoint from the simplex $\sigma$.  The matrix $P$ of WP pairings for a combined short and relative length basis  $\{\lambda_{\alpha},J\lambda_{\alpha},\grad\ell_{\beta}\}_{\alpha\in\sigma,\,\beta\in\tau}$ determines a germ at $p$ of a Lipschitz map  $P(q)=P(p)+O(d_{WP}(q,p))$ from $\Tbar$ into a real general linear group $GL(\mathbb R)$.  
\end{lemma}

We now consider the geometry of the frames provided by combined short and relative length bases.  Lemma \ref{paircont} provides the structure for a (Lipschitz) extension of the tangent bundle of $\caT$ over the stratum $\caT(\sigma)$ of (formal) codimension $2|\sigma|$.  The extension is a bundle over $\Tbar$ of rank $\dim_{\mathbb R}\caT$.  We note from Theorem \ref{gradpair} and formula (\ref{cos}) that the projections of the gradients $\{\grad\ell_{\beta}\}_{\beta\in\tau}$ onto the span of $\{\lambda_{\alpha},J\lambda_{\alpha}\}_{\alpha\in\sigma}$ have length bounded as $O(\sum_{\alpha\in\sigma}\ell_{\alpha}^{3/2})$.  The span of the first set of vectors is almost the orthogonal complement of the second set of vectors.  We now show with the following considerations that the WP connection almost splits near $\caT(\sigma)$ for the pair of subspaces.

We consider $L^{\infty}$ supremum norm bounds for Beltrami differentials as an ingredient to bound the connection.  For a geodesic $\alpha$ we conjugate the group $\Gamma$ for the geodesic to correspond to the deck transformation $A:z\rightarrow e^{\ell_{\alpha}}z$ with imaginary axis $\tilde\alpha$ and consider the coset sums  (\ref{series}).  The transformation $A$ has fundamental domain $\mathcal F=\{1\le|z|<e^{\ell_{\alpha}}\}$ with general collar $\mathcal C_a=\{1\le|z|<e^{\ell_{\alpha}},a\le\arg z\le\pi -a\}$ and collar exterior $\mathcal E_a=\mathcal F-\mathcal C_a$.  The associated distance from $\arg z=a$ to the axis $\tilde\alpha$ of $A$ is $\log \csc a+|\cot a|$.  The collar Theorem \ref{collars} provides for $\ell_{\alpha}$ small that the collar of approximate width $\log 4/\ell_{\alpha}$ embeds into the quotient $\mathbb H/\Gamma$; the collar $\mathcal C_{\ell_{\alpha}/2}$ embeds into the quotient.

\begin{lemma}
\label{gradbd}
The Beltrami differentials for gradients of simple geodesic-lengths are bounded as follows.  For a short geodesic $\alpha$ with lift the imaginary axis the supremum of $\grad\ell_{\alpha}-\frac{2}{\pi}\overline{(\frac{dz}{z})^2}(ds^2)^{-1}$ on $\mathcal C_{\ell_{\alpha}/2}$ and of $\grad\ell_{\alpha}$ on the complement of $\mathcal C_{\ell_{\alpha}/2}$ in $\mathbb H/\Gamma$ are bounded as $O(\ell_{\alpha}^2)$.   For a simple geodesic $\beta$ disjoint from the short geodesics the gradient $\grad\ell_{\beta}$ is supremum bounded as $O(\ell_{\beta})$.  The remainder term constants depend only on the constant for short geodesic length.  
\end{lemma} 
\prf We follow the approach of \cite[Lemma 3.13]{Wlbhv} and \cite{Wlspeclim}. We bound the coset sums in $L^1$ and then combine the mean value estimate with the relation between injectivity radius and exponential distance for a pointwise bound.

The main matter is to understand the contribution of a collar complement in $\mathbb H$.  We consider the collar exterior integral $\int_{\mathcal E_a}|\frac{dz}{z}|^2=2a\ell_{\alpha}$ and use the formula for distance to $\tilde\alpha$ to give the  bound $ce^{-d(\tilde\alpha,\partial\mathcal E_a)}\ell_{\alpha}$ for the integral for $a$ small.   We then apply the mean value estimate for the sum $\mathcal S=\sum_{\{B\}}\overline{B^*(\frac{dz}{z})^2}(ds^2)^{-1}$ for summand points $\{B(p)\}$ contained in the collar exterior $\mathcal E_a$ to find the bound $|\mathcal S(p)|\le c\, inj(p)^{-1}e^{-d(\tilde\alpha,\partial\mathcal E_a)}\ell_{\alpha}$. Lemma \ref{enhcollar} will provide a bound for the injectivity radius.  The estimate for the gradient is reduced to considering $d(\tilde\alpha,\{B(p)\})$.  

We first consider a short geodesic by taking $\{B\}=\langle A\rangle\backslash\Gamma-\langle A\rangle$ and a point $p\in\mathcal C_{\ell_{\alpha}/2}$.  The points $\{B(p)\}$ are contained in $\mathcal E_{\ell_{\alpha}/2}$ and a curve from a $B(p)$ to the imaginary axis $\tilde\alpha$ projects to $\mathbb H/\Gamma$ to a curve that at least leaves the collar and re enters the collar to connect to the core geodesic.  The distance $d(B(p),\tilde\alpha)$ is at least the sum of the distance to the collar boundary and the collar width $\log 4/\ell_{\alpha}$.  The resulting bounds provide the first $O(\ell_{\alpha}^2)$ estimate.  We next consider $p\in\mathcal E_{\ell_{\alpha}/2}$ with $inj(p)\ge c_0>0$.  The injectivity radius bound provides that $\{B(p)\}$ is contained within a fixed neighborhood of $\mathcal E_{\ell_{\alpha}/2}$.  The resulting distance to $\tilde\alpha$ is at least the collar width.  The bounds provide the second $O(\ell_{\alpha}^2)$ estimate.  We finally consider that $p\in\mathcal E_{\ell_{\alpha}/2}$ and is contained in a cusp region or in the collar for a geodesic distinct from $\alpha$.  The distance $d(B(p),\tilde\alpha)$ is at least the sum of the distance to the cusp or collar boundary and the collar width.  The bounds provide the third $O(\ell_{\alpha}^2)$ estimate.  The bounds for a short geodesic are complete.  

We consider a simple geodesic $\beta$ disjoint from the short geodesics.  By Lemma \ref{separ} the geodesic $\beta$ is disjoint from the cusp regions and collars about the short geodesics. The distance from a point to $\beta$ is at least the distance to the cusp region or collar boundary.  The bounds provide the desired $O(\ell_{\beta})$ estimate.  The proof is complete.

We are ready to consider continuity of the Hessian of a geodesic-length disjoint from the short geodesics.  We establish continuity by estimating in terms of the lengths of short geodesics.  Consider quantities in a neighborhood of a point of $\caT(\sigma)\subset\Tbar$.  We introduce the convention that on $\caT(\sigma)$ the Hessian $\hess \ell_{\beta}(\grad\ell_{\gamma},\grad\ell_{\eta})$ vanishes for geodesics $\beta,\gamma$ and $\eta$ not on a common component (part) of the limit.

\begin{lemma}
\label{bhesscon}
The quantity $\hess \ell_{\beta}(\grad\ell_{\gamma},\grad\ell_{\eta})$ is continuous in a neighborhood of a point of $\mathcal T(\sigma)$ for geodesics $\beta,\gamma$ and $\eta$ disjoint from the short geodesics.  
\end{lemma}
\prf  We consider a sequence in $\Tbar$ converging to a point of $\caT(\sigma)$.  Conjugate the $PSL(2;\mathbb R)$ representations to arrange that each $\beta$ has the imaginary axis $\tilde\beta$ as a lift.  We suppress subscripts and write $\Gamma$ for a group in the sequence.  For $\langle B\rangle\subset\Gamma$, the stabilizer of $\tilde\beta$, the corresponding fundamental domains $\mathcal F$ converge. The formula $\frac{4}{\pi}\Re\int_{\caF}\mu\dot{g}_{\zeta}\,dE$ for the Hessian of $\ell_{\beta}$ is presented in Lemma \ref{ldot}.  By Lemma \ref{gradbd} the Beltrami differentials for $\grad\ell_{\gamma}$ and $\grad\ell_{\eta}$ are uniformly bounded in $L^{\infty}$.  Provided the Beltrami differentials converge pointwise then standard estimates provide that the $\zeta$-derivatives of the deformation vector fields $\dot{g}$ are bounded and converge in $L^1(\caF)$, \cite[Sec. 5.3, Coro.]{AB}.  In summary pointwise convergence of the gradients $\grad\ell_{\gamma}$ and $\grad\ell_{\eta}$ of geodesic-lengths disjoint from the short geodesics will provide the desired continuity of the Hessian.  

We establish the pointwise convergence.  We conjugate the $PSL(2;\mathbb R)$ representations a second time to arrange that the $\gamma$ have the imaginary axis as a lift.  The Chabauty convergence of representations and from Lemma \ref{enhmv} the $L^1$ collar exterior bound $\int_{\mathcal E_a}|\frac{dz}{z}|^2=2a\ell_{\gamma}$ for the partial sum $\sum_{\{B\}} \overline{B^*(\frac{dz}{z})^2}(ds^2)^{-1}$ with $\{B(p)\}\subset\mathcal E_a$ provides that the $\grad\ell_{\eta}$ converge pointwise.  For $\beta,\gamma$ and $\eta$ on a common component of the Chabauty limit the conjugations differ by a convergent sequence in $PSL(2;\mathbb R)$.  For $\beta,\gamma$ and $\eta$ on a common component the gradients converge pointwise.  Next for $\beta,\gamma$ and $\eta$ not on a common component of the Chabauty limit let $\gamma$ be the geodesic not on the component of $\beta$.  A curve from $\beta$ to $\gamma$ necessarily crosses the collar of a short geodesic.  The distance from $\beta$ to $\gamma$ is at least $\min_{\alpha\in\sigma}2\log 1/\ell_{\alpha}$.  Equivalently in the upper half plane the total lift of $\beta$ is contained in a collar exterior $\mathcal E_{c\ell_{\alpha}^2}$ for $\gamma$.  In particular on  fixed-distance neighborhoods of $\beta$ the Beltrami differentials for the $\grad\ell_{\gamma}$ are bounded in $L^1$ as $O(\ell_{\alpha}^2)$.  From Lemma \ref{enhmv} the $\grad\ell_{\gamma}$ converge pointwise to zero.  Since the deformation vector fields $\dot{g}_{\zeta}$ are bounded in $L^1$, convergence of the $\ddot{\ell}_{\beta}$-integral to zero now follows.  The proof is complete. 

We bound the Hessian evaluated on the gradient of short geodesic-lengths.  For the purpose of bounding $\mathbb P_{\beta}$ we introduce eigenfunctions of the Laplace-Beltrami operator.  Considerations are based on two properties of eigenfunctions.  First, given $\epsilon$ positive there is a constant $c(\epsilon)$ such that for a $2$-eigenfunction $f$ on a domain in $\mathbb H$ containing a metric ball $B(p;\epsilon)$ the mean value property \cite[Coro. 1.3]{Fay} is satisfied
\[
f(p)=c(\epsilon)\int_{B(p;\epsilon)}f\,dA.
\]
The constant is evaluated by considering an initially positive radial eigenfunction $k$.  The equation $Dk=2k$ provides for an initial strict local minimum.  The local minimum and equation 
$k(p)=c(\epsilon)\int_{B(p;\epsilon)}k\,dA$ provide the inequality $0<Area(B(p;\epsilon))\,c(\epsilon)<1$.  Second, for $z\in\mathbb H$ with polar form $z=re^{i\theta}$ the increasing non negative function $u(\theta)=1-\theta\cot \theta$ is a $2$-eigenfunction.  The modification $v(\theta)=\min\{{u(\theta),u(\pi-\theta)}\}$ is a continuous function and a $2$-eigenfunction on the complement of the imaginary axis.  The initial expansion $u(\theta)=\frac13\theta^2+O(\theta^4)$ provides that $v(\theta)$ is comparable to $\sin^2\theta$.  In particular for the geodesic $\beta$ corresponding to a deck transformation $B$ with imaginary axis $\tilde\beta$ the coset sum $\mathbb M_{\beta}(p)=\sum_{C\in\langle B\rangle\backslash\Gamma}v(\theta(C(p))$ is comparable to $\mathbb P_{\beta}(p)$.  The $\Gamma$-invariant function $\mathbb M_{\beta}$ satisfies the mean value property on 
$\mathbb H-\cup_{C\in\langle B\rangle\backslash\Gamma}C^{-1}(\tilde\beta)$ (for $\epsilon$-balls contained in the region.)

We consider $\mathbb P_{\beta}$ on collars for short geodesics.  A uniform bound follows from the mean value estimate Lemma \ref{enhmv}. To further bound $\mathbb P_{\beta}$ or equivalently $\mathbb M_{\beta}$ in a collar represent the collar in $\mathbb H$ with coordinate $z=re^{i\theta}$ as $\mathcal C_{\ell_{\alpha}}=\{1\le r<e^{\ell_{\alpha}}, \ell_{\alpha}\le \theta\le \pi-\ell_{\alpha}\}$.  The uniform bound for $\mathbb M_{\beta}$ and expansion $u=\pi/\theta$ for $\theta$ close to $\pi$ combine to provide a positive constant $c$ such that $h=c\,\ell_{\alpha}(u(\theta)+u(\pi-\theta))-\mathbb M_{\beta}$ is positive on the boundaries $\theta=\ell_{\alpha},\pi-\ell_{\alpha}$ of the collar.

We consider the behavior of the $\langle B\rangle$-invariant function $h$ in the collar.  For all small $\epsilon$ an $\epsilon$-neighborhood of the collar is disjoint from the $\Gamma$-orbit of $\tilde\beta$.  The function $h$ satisfies the mean value property for points of the collar.  For a small $\epsilon$ the inequality $Area(B(p;\epsilon))\,c(\epsilon)<1$ excludes the possibility of $\max_{\mathcal C_{\ell_{\alpha}}}h$ being realized at points $p$ for which $B(p;\epsilon)\subset\mathcal C_{\ell_{\alpha}}$.  A maximum is realized only on the boundary.  In conclusion $c\,\ell_{\alpha}(u(\theta)+u(\pi-\theta))$ provides a uniform upper bound for $\mathbb P_{\beta}$ on the collar. 

We are ready to consider the Hessian evaluated on the gradients of short geodesic-lengths.
\begin{lemma}
\label{babd}
The quantities $\hess\ell_{\beta}(V,V)$  for $V=\grad\ell_{\alpha}$ and $J\grad\ell_{\alpha}$ are bounded as $O(\ell_{\alpha}^2)$ with a remainder constant depending only on the bounds for the length of $\alpha$ and $\beta$.
\end{lemma}
\prf The Hessian by Theorem \ref{hessbd} is bounded by the WP pairing with the weight function $\mathbb P_{\beta}$.  From the Cauchy Schwartz inequality and Lemma \ref{gradbd} it is enough to separately consider the contribution of the principal term $\frac{2}{\pi}\overline{(\frac{dz}{z})^2}(ds^2)^{-1}$ on the collar $\mathcal C_{\ell_{\alpha}}$ and a remainder term bounded in $L^{\infty}$ by $O(\ell_{\alpha}^2)$.  In fact since $\mathbb P_{\beta}$ has bounded $L^1$ norm the contribution of the remainder term to $\hess\ell_{\beta}$ is $O(\ell_{\alpha}^4)$, a suitable bound.

We consider the principal term for the collar.  The term for the Beltrami differential is $\frac{2}{\pi}\overline{(\frac{dz}{z})^2}(ds^2)^{-1}$.  A bound for the weight function is $c\,\ell u(\theta)$ and the resulting comparison integral is
\[
\int^{e^{\ell}}_1\int^{\pi-\ell}_{\ell}r^2\sin^2\theta\,\ell(1-\theta\cot\theta)\frac{|dz|^2}{|z|^4}
= \ell^2 \int^{\pi-\ell}_{\ell} \sin^2\theta(1-\theta\cot\theta)\,d\theta=O(\ell^2).
\]  The proof is complete.

We are ready to describe the WP connection in a neighborhood of a point of $\caT(\sigma)\subset\Tbar$ with the short geodesic-lengths $\sigma$ and relative length basis $\tau$.  The WP metric is K\"{a}hler and consequently $J$ is parallel.  
\begin{theorem}
\label{WPconn}
The WP connection for a combined frame of gradients $\{\lambda_{\alpha},J\lambda_{\alpha},\grad\ell_{\beta}\}_{\alpha\in\sigma,\beta\in\tau}$ for $\alpha,\alpha'\in\sigma$, $\beta,\beta'\in\tau$  satisfies: \\
$D_{J\lambda_{\alpha}}\lambda_{\alpha}=3/(2\pi\ell_{\alpha}^{1/2})J\lambda_{\alpha}+O(\ell_{\alpha}^{3/2}),$\;
$D_{\lambda_{\alpha'}}\lambda_{\alpha}=O(\ell_{\alpha}^{3/2})$,\;\\   
$D_{J\lambda_{\alpha'}}\lambda_{\alpha}=O(\ell_{\alpha}(\ell_{\alpha'}^{3/2}+\ell_{\alpha}^{1/2}))$, for 
$\alpha\ne\alpha'$,\;  
$D_{\grad\ell_{\beta}}\lambda_{\alpha}=O(\ell_{\alpha})$, \\ 
$D_{\lambda_{\alpha}}\grad\ell_{\beta}=O(\ell_{\alpha}^{1/2})$,\; $D_{J\lambda_{\alpha}}\grad\ell_{\beta}= 
O(\ell_{\alpha}^{1/2})$,\; and \; $D_{\grad\ell_{\beta'}}\grad\ell_{\beta}$ 
is continuous at $\caT(\sigma)$  with value zero if $\beta$ and $\beta'$ lie on distinct parts.
The remainder term constants depend only on the constant for short geodesic-length.
\end{theorem}
\prf  The expansions follow from combining the above considerations.  The proof is complete.

\section{Approximating WP geodesics}
\label{WPgeoapp}
We show in general that a geodesic in a stratum $\mathcal T(\sigma)\subset\Tbar$ is appropriately $C^1$-approximated 
by geodesics in the Teichm\"{u}ller space $\mathcal T$.  In Section \ref{horocon} we use the property to approximate geodesics for punctured surfaces by geodesics for compact surfaces, and to then analyze the variation of the distance between horocycles.  To understand the geodesic approximation near a point of $\mathcal T(\sigma)$ we write the geodesic equation in terms of the frame of a combined short and relative length basis and apply the estimates of Sections \ref{gradhess} and \ref{conWP}.  The approximation is presented in Theorem \ref{approx}.    

We describe the setup.  The combined short and relative length basis $\{\lambda_{\alpha},J\lambda_{\alpha},\grad\ell_{\beta}\}_{\alpha\in\sigma,\beta\in\tau}$ provides a frame $\mathcal F$ for the tangent space $T\mathcal T$ near a point of $\mathcal T(\sigma)$, and the relative length basis $\{\grad\ell_{\beta}\}_{\beta\in\tau}$ provides a frame for the tangent space $T\mathcal T(\sigma)$ near the point.  We denote the general element of the frame $\mathcal F$ as $v_{\kappa}$ and also write $\{v_{\kappa}\}$ for the frame.  Considerations will involve the pairs of gradients $\lambda_{\alpha},J\lambda_{\alpha}$ for $\alpha\in\sigma$.  We write $(J)\lambda_{\alpha}$ or $v_{(J)\alpha}$ to denote the pair.  We begin considerations by expressing the tangent field of a geodesic $\gamma(t)$ in $\mathcal T$ or $\mathcal T(\sigma)$ in terms of the frame $\{v_{\kappa}\}$ as follows $\dot\gamma=(v_{\kappa})^t(c_{\kappa})$ for a column vector of coefficients and the frame elements (for $\gamma\subset\mathcal T(\sigma)$ the subframe $\{v_{\beta}\}_{\beta\in\tau}$ locally spans $T\mathcal T(\sigma)$ and the coefficients $c_{(J)\alpha}$ are set equal to zero.)  We now have the equation
\begin{equation}
\label{bdrv1}
(\dot\ell_{\beta}(\gamma))=(\langle v_{\beta},\dot\gamma\rangle)=(\langle v_{\beta},v_{\kappa}\rangle)(c_{\kappa})\quad \mbox{for}\ \beta\in\tau
\end{equation}
for the column vectors of derivatives and the matrix of gradient pairings.  We have from the geodesic equation 
$D_{\dot\gamma}\dot\gamma=0$ and differentiating along the geodesic the equation
\begin{equation}
\label{bdrv2}
(\ddot\ell_{\beta}(\gamma))=(\dot\gamma\langle v_{\beta},\dot\gamma\rangle)=(\langle D_{\dot\gamma}v_{\beta},\dot\gamma\rangle)= (c_{\rho})^t(\langle D_{v_{\rho}}v_{\beta},v_{\kappa}\rangle)(c_{\kappa})\quad\mbox{for}
\ \beta\in\tau
\end{equation}
for the appropriate column vectors.

To relate geodesics in $\mathcal T$ to geodesics in $\mathcal T(\sigma)$ we introduce a coordinate projection and a coefficient projection for vectors given as a sum of frame elements.  We introduce coordinates for a neighborhood of $p\in\mathcal T(\sigma)$ as follows.  Begin with a partition $\mathcal P$ for the reference surface $F$ with $\sigma\subset\mathcal P$.  Short geodesic-lengths and Fenchel-Nielsen angles in combination with the relative length basis lengths $\{\ell_{\alpha},\vartheta_{\alpha},\ell_{\beta}\}_{\alpha\in\sigma,\beta\in\tau}$ provide coordinates for a neighborhood $\mathcal U$ of $p$.  We define locally the {\em coordinate projection} by associating to a point $q\in\mathcal U$ with coordinates $\{\ell_{\alpha}(q),\vartheta_{\alpha}(q),\ell_{\beta}(q)\}$ the projected point $q_{\dagger}\in\mathcal T(\sigma)$ with coordinates $\{\ell_{\alpha}=0, \ell_{\beta}(q)\}$.  We define locally the 
{\em coefficient projection} by associating to the element of the span of $\mathcal F$ with coefficients $(a_{\kappa})_{\kappa\in\sigma\cup J\sigma\cup\tau}$ the element with coefficients $(a_{\kappa})_{\kappa\in\tau}$ (the coefficients with indices in $\sigma\cup J\sigma$ are set equal to zero.)  We write $(a_{\kappa})_{\dagger}$ for the projection of an element $(a_{\kappa})$ of the span of $\mathcal F$ and write $(\langle v_{\eta},v_{\kappa}\rangle)_{\dagger}$ for the projection $(\langle v_{\eta},v_{\kappa}\rangle)_{\eta,\kappa\in\tau}$ of the pairing matrix.  We now use the projections to relate geodesics.  A unit-speed geodesic $\gamma(t)$ in $\mathcal T$ close to $\mathcal T(\sigma)$ has tangent field described in terms of the frame $\mathcal F$ and projections as follows.
\begin{lemma}
\label{geodproj}
For a point in $\mathcal T(\sigma)$ there is a neighborhood $\mathcal U$ in $\Tbar$ and a positive value $t_0$ such that a unit-speed geodesic $\gamma(t)$ with initial point in $\mathcal U$ and $\langle\lambda_{\alpha},\dot\gamma(0)\rangle=\langle J\lambda_{\alpha},\dot\gamma(0)\rangle=0$, $\alpha\in\sigma$, is described as follows.  The root geodesic-length functions $\ell_{\alpha}^{1/2}$ on the geodesic satisfy for $0\le t\le t_0$ the inequalities $\ell_{\alpha}^{1/2}(0)\le\ell_{\alpha}^{1/2}(t)\le2\ell_{\alpha}^{1/2}(0)$ with the quantities 
$\langle\lambda_{\alpha},\dot\gamma\rangle,\ \langle J\lambda_{\alpha},\dot\gamma\rangle$ bounded as $O(\ell_{\alpha}^{3/4})$.  The derivatives of the geodesic-length functions $\ell_{\beta}$ along the geodesic satisfy
\begin{align*}
(\dot\ell_{\beta}(\gamma))=(\langle v_{\beta},v_{\kappa}\rangle)_{\dagger}(c_{\kappa})_{\dagger}\ +\ O(\sum_{\alpha\in\sigma}\ell_{\alpha}) \quad\quad\mbox{and} \\
(\ddot\ell_{\beta}(\gamma))=(c_{\rho})_{\dagger}^t(\langle D_{v_{\rho}}v_{\beta},v_{\kappa}\rangle)_{\dagger}(c_{\kappa})_{\dagger}\ +\ O(\sum_{\alpha\in\sigma}\ell_{\alpha}) \notag
\end{align*}
for the tangent field expressed as $\dot\gamma=(v_{\kappa})^t(c_{\kappa})$ in terms of the frame $\mathcal F$.  The value $t_0$ and remainder term constants depend only on the choice of the set $\mathcal U$.
\end{lemma}
\prf The considerations are a collection of estimates.  We first consider $\langle(J)\lambda_{\alpha},\dot\gamma\rangle$.  From the initial vanishing of the derivative and the convexity of root geodesic-length functions it follows that $\ell_{\alpha}^{1/2}$ is increasing on $\gamma$ for $t$ positive.  On an interval of $t$ positive the functions $\ell_{\alpha}^{1/2}$ are bounded by their endpoint values.  We introduce the function $F(t)=\langle\lambda_{\alpha},\dot\gamma(t)\rangle^2+\langle J\lambda_{\alpha},\dot\gamma(t)\rangle^2$ and use Theorem \ref{WPconn} to find an expansion for its derivative.  The principal terms combine to cancel and it follows that $\frac{dF}{dt}$ is $O(\ell_{\alpha}^{3/2})$ with a uniform constant provided $\ell_{\alpha}\le c_0$.  Since $F(0)=0$ it follows that $F(t)$ is $O(t\ell_{\alpha}^{3/2})$ with a uniform constant provided $\ell\le c_0$.  Further since $\langle\lambda_{\alpha},\dot\gamma(t)\rangle^2\le F(t)$ it follows that $\ell_{\alpha}^{1/2}(t) \le \ell_{\alpha}^{1/2}(0)+O(t^{3/2}\ell_{\alpha}^{3/4}) \le \ell_{\alpha}^{1/2}(0) +O(t^{3/2})$ with a uniform constant on any interval with $\ell_{\alpha}\le c_0$.  It follows for $\ell_{\alpha}(0)$ sufficiently small that the interval with $\ell_{\alpha}\le 1$ contains a fixed initial open interval in $t$.  (The values $\ell_{\alpha}$, $1/\ell_{\alpha}$ are bounded on $\gamma$.  We apply Lemma \ref{paircont} to bound the gradients $\grad \ell_{\beta}$ on $\gamma$, and to bound the values $\ell_{\beta}$, $1/\ell_{\beta}$.  It follows that the geodesic $\gamma$ exists for a uniform interval in $t$.)  It further follows that there is an initial subinterval on which the remainder term $O(t^{3/2}\ell_{\alpha}^{3/4})$ is explicitly bounded by $\ell_{\alpha}^{1/2}/2$ which results in the inequality $\ell_{\alpha}^{1/2}(t)\le2\ell_{\alpha}^{1/2}(0)$ for the subinterval.  The estimates for $\ell_{\alpha}^{1/2}$ and $\langle (J)\lambda_{\alpha},\dot\gamma\rangle$ are complete.

We proceed and bound the coefficients $c_{\alpha}$ and $c_{J\alpha}$ as $O(\ell_{\alpha}^{3/4})$.  The tangent field $\dot\gamma$ has unit length.  A uniform bound for the coefficients follows from the defining equation 
$\dot\gamma=(v_{\kappa})^t(c_{\kappa})$ and Lemma \ref{paircont}.  Next from Theorem \ref{gradpair} the off-diagonal pairings $\langle v_{\alpha},v_{\kappa}\rangle$, $\kappa\ne\alpha$, and $\langle v_{J\alpha},v_{\kappa}\rangle$, $\kappa\ne J\alpha$, are bounded at least as $O(\ell_{\alpha}^{3/2})$.  We combine estimates to find that 
$\langle v_{\alpha},\dot\gamma\rangle=\langle v_{\alpha},v_{\alpha}\rangle c_{\alpha}+ O(\ell_{\alpha}^{3/2})$  
and 
$\langle v_{J\alpha},\dot\gamma\rangle=\langle v_{J\alpha},v_{J\alpha}\rangle c_{J\alpha}+ O(\ell_{\alpha}^{3/2})$ with uniform remainder constants.  The desired bounds for $c_{\alpha}$ and $c_{J\alpha}$ now follow from the bounds for $\langle (J)\lambda_{\alpha},\dot\gamma\rangle$ and Theorem \ref{gradpair}.  The coefficient estimates are complete.       
  
We proceed and bound the derivatives of $\ell_{\beta}$ along $\gamma$.  We consider (\ref{bdrv1}) and gather terms with $\kappa$ in $\sigma\cup J\sigma$ to form a remainder term.  Theorem \ref{gradpair} provides for the first expansion with a uniform remainder constant.  We consider (\ref{bdrv2}) and gather terms with at least one of $\rho$ or $\kappa$ in $\sigma\cup J\sigma$ to form a remainder term.  Lemma \ref{babd} and the estimates for $c_{(J)\alpha}$ provide for the second expansion with a uniform remainder constant.  The proof is complete.

We consider with Lemma \ref{geodproj} a unit-speed geodesic $\gamma_0(t)$ contained in $\mathcal T(\sigma)$ with $\mathcal U$ a suitable neighborhood in $\Tbar$ of the initial point and $t_0$ a suitable value.  Given $\delta$ positive we consider an approximating geodesic $\gamma(t)$ contained in $\mathcal T$ with 
initial point in $\mathcal U$, with $\langle\lambda_{\alpha},\dot\gamma(0)\rangle=\langle J\lambda_{\alpha},\dot\gamma(0)\rangle=0$, for $\alpha\in\sigma$, and with $|\dot\ell_{\beta}(\gamma(0))-\dot\ell_{\beta}(\gamma_0(0))|<\delta$, for $\beta\in\tau$.  In general we write $(\ell_{\beta}(\gamma),\dot\ell_{\beta}(\gamma))$ and $(\ell_{\beta}(\gamma_0),\dot\ell_{\beta}(\gamma_0))$ to denote the column vectors of values and first derivatives along the geodesics.  
We combine the above considerations to derive integral-equations for the vectors $(\ell_{\beta}(\gamma),\dot\ell_{\beta}(\gamma))$ and $(\ell_{\beta}(\gamma_0),\dot\ell_{\beta}(\gamma_0))$ for the relative length basis.  We consider preliminaries.  The coordinate projection satisfies $\ell_{\beta}(\gamma_{\dagger}(t))=\ell_{\beta}(\gamma(t))$ for $\beta\in\tau$ and consequently differentiation in $t$ commutes with the coordinate projection.  We write for the relative length basis
\[
P_{\dagger}=(\langle v_{\eta},v_{\kappa}\rangle)_{\dagger}\quad\mbox{and}\quad A=P_{\dagger}^{-1}(\langle D_{v_{\rho}}v_{\eta},v_{\kappa}\rangle)_{\dagger}P_{\dagger}^{-1}
\]
for the symmetric matrix of gradient pairings and the conjugated matrix of covariant derivatives.  By Lemma \ref{bhesscon} the covariant derivatives $\langle D_{v_{\rho}}v_{\eta},v_{\kappa}\rangle$ for $\rho,\eta,\kappa\in\tau$ are continuous in a neighborhood of the stratum $\mathcal T(\sigma)$.  We can further select the neighborhood $\mathcal U$ to arrange that the intersection $\mathcal U\cap \mathcal T(\sigma)$ is relatively compact in $\mathcal T(\sigma)$.  It then follows from Lemma \ref{bhesscon} that $A(\gamma)=A(\gamma_{\dagger})+o(1)$ with a uniform remainder, in particular the matrix $A$ evaluated on $\gamma$ and on the coordinate projection $\gamma_{\dagger}$ differ by a quantity uniformly bounded along $\gamma$ and $\gamma_{\dagger}$.  We conclude from Lemmas \ref{paircont} and \ref{geodproj} for a suitable neighborhood $\mathcal U$ that for $0\le t\le t_0$
\[
(\ddot\ell_{\beta}(\gamma))=(\dot\ell_{\rho}(\gamma))^tA(\gamma_{\dagger})(\dot\ell_{\kappa}(\gamma))\ +\ o(1)
\]
for $\rho,\kappa\in\tau$ and on integration conclude that
\begin{equation}
\label{gameq}
(\dot\ell_{\beta}(\gamma(t)))=\bigl(\dot\ell_{\beta}(\gamma(0))\ +\ \int_0^t\bigl((\dot\ell_{\rho}(\gamma))^t A(\gamma_{\dagger})(\dot\ell_{\kappa}(\gamma))\bigr)(s)\,ds\ +\  o(1)\bigr).
\end{equation}
The geodesic $\gamma_0$ is contained in the stratum $\mathcal T(\sigma)$ with the local coordinates $(\ell_{\beta})$ for which an exact relation is satisfied for $0\le t\le t_0$
\[
(\ddot\ell_{\beta}(\gamma_0))=(\dot\ell_{\rho}(\gamma_0))^t A(\gamma_0)(\dot\ell_{\kappa}(\gamma_0))
\]
for $\rho,\kappa\in\tau$ and on integration conclude that
\begin{equation}
\label{gam0eq}
(\dot\ell_{\beta}(\gamma_0(t)))=\bigl(\dot\ell_{\beta}(\gamma_0(0))\ +\ \int_0^t\bigl((\dot\ell_{\rho}(\gamma_0))^t A(\gamma_0 )(\dot\ell_{\kappa}(\gamma_0))\bigr)(s)\,ds\bigr) .
\end{equation}
We are ready to describe the $C^1$-approximation of the geodesic $\gamma_0$ by the geodesic $\gamma$ in terms of the short geodesics $\sigma$ and the relative length basis $\tau$.
\begin{theorem}
\label{approx}
For a point $p$ in $\mathcal T(\sigma)$ and $\epsilon$ positive there exist a neighborhood $\mathcal U$ in $\Tbar$, and positive values $\delta,t_0$ as follows.  For a unit-speed geodesic $\gamma_0$ contained in $\mathcal T(\sigma)$ with $\gamma_0(0)=p$ then a geodesic $\gamma$ with $\gamma(0)$ in $\mathcal U$, with $\langle\lambda_{\alpha},\dot\gamma(0)\rangle=\langle J\lambda_{\alpha},\dot\gamma(0)\rangle=0$, $\alpha\in\sigma$, and with 
$\|(\dot\ell_{\beta}(\gamma(0))-\dot\ell_{\beta}(\gamma_0(0)))\|<\delta$ approximates as follows.   
The geodesics $\gamma,\gamma_0$ exist on the interval $-t_0\le t\le t_0$.  The root geodesic-length functions $\ell_{\alpha}^{1/2}$, $\alpha\in\sigma$, on $\gamma(t)$ for $-t_0\le t\le t_0$ satisfy the inequalities $\ell_{\alpha}^{1/2}(0)\le\ell_{\alpha}^{1/2}(t)\le2\ell_{\alpha}^{1/2}(0)$ with the quantities 
$\langle\lambda_{\alpha},\dot\gamma\rangle,\ \langle J\lambda_{\alpha},\dot\gamma\rangle$ bounded as $O(\ell_{\alpha}^{3/4})$. The relative length basis functions for $-t_0\le t\le t_0$ satisfy the inequality 
$\|(\ell_{\beta}(\gamma),\dot\ell_{\beta}(\gamma))- (\ell_{\beta}(\gamma_0),\dot\ell_{\beta}(\gamma_0))\|<\epsilon$.  The remainder term constant depends only on the choice of $p$ and $\epsilon$.  
\end{theorem}
\prf  The bounds for the short geodesic-lengths and the existence interval for the geodesics are provided in Lemma \ref{geodproj}.  We now refer to the integral equations (\ref{gameq}) and (\ref{gam0eq}) to apply the basic inequality for first-order differential equations.   Observations are in order.  For the integral equations (\ref{gameq}) and (\ref{gam0eq}) the matrix-valued function $A$ is evaluated on $\mathcal U\cap\mathcal T(\sigma)$, since the coordinate projection $\gamma_{\dagger}$ lies in  and $\gamma_0$ is given in $\mathcal T(\sigma)$.  The WP metric on $\mathcal T(\sigma)$ is smooth and the intersection $\mathcal U\cap\mathcal T(\sigma)$ is given above as relatively compact.  The function $f(\ell_{\beta},\dot\ell_{\beta})= (\dot\ell_{\beta})^tA(\dot\ell_{\beta})$ of the values $(\ell_{\beta},\dot\ell_{\beta})$ is uniformly Lipschitz 
on $\mathcal U\cap\mathcal T(\sigma)$.  We can apply Theorem 2.1 of \cite{CdLv} to obtain the desired bounds for the relative length basis.  The proof is complete.

\section{Convexity of the distance between horocycles}
\label{horocon}
We show for a pair of cusps on a hyperbolic Riemann surface that the distance between horocycles is a WP strictly convex function on the Teichm\"{u}ller space.  The convexity is presented as an application of the convexity of geodesic-length functions \cite{Wlnielsen, Wlbhv}, and the approximation of geodesics in a stratum described in Theorem \ref{approx}.  To this purpose for a hyperbolic geodesic $\omega_0$ between cusps $p,q$ for a surface $R_0$ we introduce the double across the cusps $R_0\cup_{p,q}\overline{R_0}$ and open the cusps $p,q$ (considered as nodes) to obtain an approximating surface $R$ with an approximating closed geodesic $\omega$.  The closed geodesic $\omega$ consists of segments in the thick regions of $R$ and segments crossing the collars $\mathcal C(p)$ and $\mathcal C(q)$. The length $\ell_{\omega}-2w(p)-2w(q)$ of $\omega$ with twice the collar widths subtracted approximates twice the distance $2\tilde\ell_{\omega_0}$ between horocycles.  We use the description of distance in terms of closed geodesics and collar widths to establish the convexity in Theorem \ref{horowpcon}.  

We describe the setup.  Consider a surface $R_0$ with a collection of cusps $\mathcal N$.  For the conjugate complex structure $\overline{R_0}$ consider the union $\DR$ as a Riemann surface with nodes where the elements of $\mathcal N$ on $R_0$ and $\overline{R_0}$ are formally paired to form nodes.  Introduce the reference topological surface $F$ with the simplex $\sigma$ of the complex of curves $C(F)$ and a reference homeomorphism from $F-\cup_{\alpha\in\sigma}\alpha$ to $\DR$ with the elements of $\sigma$ corresponding to the cusp-pairs of $\mathcal N$.  The Riemann surface with nodes $\DR$ represents a point of the stratum $\mathcal T(\sigma)$ for the augmented Teichm\"{u}ller space $\Tbar(F)$ of hyperbolic structures for $F$.  

The identity map for the underlying topological surface $R_0$ induces an anti conformal reflection $\upsilon_0$ of $\DR$.  The mapping $\upsilon_0$ corresponds to an orientation reversing self homeomorphism $\upsilon$ of $F$ with 
$\upsilon\circ\upsilon=id$ and with $\upsilon$ fixing the elements of $\sigma$.  The mapping $\upsilon$ determines an involutive element of the extended mapping class group.  The mapping acts on $\mathcal T(F)$ as an anti biholomorphic WP isometry.   Following Definition \ref{relbas}, there is a local frame for the tangent bundle adapted to the 
action of $\upsilon$.  Accordingly introduce a relative length basis $\tau$ for $\DR$ with the elements of $\tau$ permuted by $\upsilon$.   The elements of $\tau$ are disjoint from the elements of $\sigma$.  The gradients $\{\lambda_{\alpha},J\lambda_{\alpha},\grad \ell_{\beta}\}_{\alpha\in\sigma,\beta\in\tau}$ provide a local frame for $T\mathcal T(F)$ in a neighborhood of $\DR$.  The mapping $\upsilon$  fixes the elements of $\sigma$ and permutes the elements of $\tau$.  The mapping acts on the local frame as follows $d\upsilon(\lambda_{\alpha})=\lambda_{\alpha}$, 
$d\upsilon(J\lambda_{\alpha})=-J\lambda_{\alpha}$ and $d\upsilon(\grad\ell_{\beta})=\grad\ell_{\upsilon(\beta)}$. 

The fixed point locus in $\Tbar(F)$ of the involution $\upsilon$ is a totally geodesic locus.  
The locus parameterizes marked structures with $\upsilon$ realized as a reflection (such marked hyperbolic structures for $F$ are in bijection to marked hyperbolic structures with geodesic boundary for the quotient $F\slash\langle\upsilon\rangle$.)  Accordingly for the WP geodesic equation, initial conditions invariant by $\upsilon$ determine a geodesic parameterizing marked structures with $\upsilon$ realized as a reflection.  

Now the considerations for Theorem \ref{horowpcon} begin with a geodesic $\gamma_{00}$ in the Teichm\"{u}ller space of $R_0$ with the geodesic containing a point representing $R_0$.  Introduce the geodesic $\gamma_0$ in the product Teichm\"{u}ller space of $\DR$ given as $\gamma_{00}\cup_{\mathcal N}\overline{\gamma_{00}}$.   For $\gamma_{00}$ describing a family of Riemann surfaces $\mathcal{R}$ with the collection of cusps $\mathcal N$ on fibers, the geodesic $\gamma_0$ describes the family of Riemann surfaces $(\mathcal{R},\overline{\mathcal{R}})$ with the cusps $\mathcal N$ on fibers of $\mathcal{R}$ and $\overline{\mathcal{R}}$ formally paired.  The geodesic $\gamma_0$ is contained in the stratum $\mathcal T(\sigma)$ and is fixed by the action of $\upsilon$.  On $\gamma_0$ the geodesic-lengths for the relative length basis $\tau$ satisfy $\ell_{\beta}(\gamma_0)=\ell_{\upsilon(\beta)}(\gamma_0)=\ell_{\beta}(\upsilon(\gamma_0))$ 
and consequently for the tangent field of the geodesic $d\ell_{\beta}(\dot\gamma_0)=d\ell_{\beta}(d\upsilon(\dot\gamma_0))$.  Now from Theorem \ref{approx} there is a neighborhood of the point $p_0=\DR$ in $\Tbar(F)$ and the geodesic $\gamma_0$ is approximated by geodesics of $\mathcal T(F)$ with 
approximating geodesics containing base points in the neighborhood.  The approximation is now specialized to specify geodesics $\gamma$ with base points $q_0$ in the neighborhood with the geodesics fixed by the action of $\upsilon$.  First, a base point $q_0$ is selected on the fixed point locus of $\upsilon$ (the locus is given by equations $\ell_{\upsilon(\beta)}=\ell_{\beta}$,  $\beta\in\tau$ and appropriate values of the Fenchel-Nielsen angles $\vartheta_{\alpha}$, $\alpha\in\sigma$.)  Second, the tangent $\mathbf t$ to $\gamma$ at $q_0$ is specified by the equations 
$d\ell_{\beta}(\mathbf t)=d\ell_{\beta}(\mathbf t_0)$, $\beta\in\tau$, for $\mathbf t_0$ the tangent to $\gamma_0$ 
at $p_0$ and the equations $\langle\lambda_{\alpha},\mathbf t\rangle =\langle J\lambda_{\alpha},\mathbf t\rangle=0$, $\alpha\in\sigma$.  The initial conditions for $\gamma$ are invariant by the action of $\upsilon$.  Third, the base point $q_0$ is close to $p_0$ provided for $\beta\in\tau$ and $\alpha\in\sigma$ the values $(\ell_{\beta}(q_0),\ell_{\alpha}(q_0))$ are close to the values $(\ell_{\beta}(p_0),0)$.  Theorem \ref{approx} now provides for approximating geodesics $\gamma$ fixed by the action of $\upsilon$ and close to $\gamma_0$.

Geodesics between cusps on the hyperbolic surface $R_0$ are approximated by closed geodesics on the approximating surfaces.  Distance between horocycles on $R_0$ is approximated by the length of geodesic segments in the complement of the $\sigma$ collars on the approximating surfaces.   The association between geodesics and segments is as follows.  A geodesic between cusps with a choice of appropriate horocycles about the cusps determines a homotopy class of an arc connecting the horocycles with endpoints relative to the horocycles.  The homotopy class does not depend on the choice of horocycles since cusp neighborhoods are fibered by horocycles.  In particular a geodesic between cusps on $R_0$ determines a geodesic and a homotopy class relative to horocycles on each component of $\DR$.  The reference homeomorphism from $F-\cup_{\alpha\in\sigma}\alpha$ to $\DR$ provides for a pair of homotopy classes relative to $\sigma$ on $F$ and so provides for a pair of homotopy classes relative to the short geodesics $\sigma$ on a marked surface $R$ of $\mathcal T(F)$ close to $\DR$.  A homotopy class relative to the short geodesics $\sigma$ on a marked surface $R$ contains a unique length minimizing representative orthogonal to $\sigma$.  For marked surfaces 
with $\upsilon$ realized as a reflection, the length minimizing representatives of the pair of homotopy classes are interchanged by the reflection and $\sigma$ is pointwise fixed by the reflection.  It follows that the length minimizing representatives combine to form a smooth closed geodesic on marked surfaces $R$ with reflection $\upsilon$ and close to $\DR$. 

The thin regions on marked surfaces $R$ consist of any cusps not in $\mathcal N$ and the collars about $\alpha$, $\alpha\in\sigma$.   The thin region of $\DR$ consists of cusps (in $\mathcal N$ and not in $\mathcal N$.)   The Chabauty topology for $\Tbar(F)$ provides that the geometry of the thick region of $R$ is close to the geometry of the thick region of $\DR$.  The closed geodesic described above consists of a segment on the thick region of $R_0$, a segment on the thick region of $\overline{R_0}$ and two segments crossing $\sigma$ collars.  The segments on the thick regions are interchanged by the reflection $\upsilon$ and the segments in the collars are stabilized by the reflection 
$\upsilon$.   The segments in the collars are orthogonal to the collar core geodesics and so are orthogonal to the collar meridians.  It follows that the length of a segment in a collar is exactly twice the width of the collar.
For $R$ with reflection $\upsilon$ and close to $\DR$ the length of a segment of the closed geodesic in the thick region is close to the length of the geodesic segment connecting horocycles on $R_0$.  In summary the length of the closed geodesic with twice the collar widths subtracted is close to twice the length of the geodesic segment connecting horocycles.

\begin{theorem}
\label{horowpcon}
For a marked surface $R_0$ with a hyperbolic geodesic $\omega_0$ between cusps and a choice of horocycle length at most $2$, the distance $\tilde\ell_{\omega_0}$ along $\omega_0$ between associated horocycles is a WP strictly convex function.  The quantity $\hess \tilde\ell_{\omega_0}$ is bounded below in terms of the WP pairing and a positive constant depending only on a positive lower bound for the geodesic-lengths for $R_0$.
\end{theorem}
\prf  As above, for $\omega_0$ a hyperbolic geodesic between not necessarily distinct cusps $p,q$ introduce the double across cusps $R_0\cup_{p,q}\overline{R_0}$.  Following Theorem \ref{approx}, consider $\gamma_0$ a geodesic containing $R_0\cup_{p,q}\overline{R_0}$ and $\gamma$ a geodesic in the corresponding Teichm\"{u}ller space $\mathcal T(F)$ with $\gamma$ close to $\gamma_0$, $\gamma$ fixed by the action of the involution $\upsilon$, and the surfaces parameterized by $\gamma$ with approximating closed geodesics $\omega$.  Following Theorem \ref{collars}, consider the collars $\mathcal C(p),\mathcal C(q)$ about the short geodesics (the opened nodes) for the surfaces parameterized by $\gamma$. The Chabauty topology provides for $\gamma$ close to $\gamma_0$ that the geometry of corresponding thick regions of surfaces are close.  In particular the function $\ell_{\omega}-2w(p)-2w(q)$ along $\gamma$ is close to the function $2\tilde\ell_{\omega_0}$ along $\gamma_0$.  From Section \ref{basics} the collar width is given as $w(\alpha)=\log 4/\ell_{\alpha}+O(\ell_{\alpha}^2)$ and the function $2\ell_{\omega}-2\log 4/\ell_{\alpha_p}-2\log 4/\ell_{\alpha_q}$ is also close to the function $2\tilde\ell_{\omega_0}$. 

We consider the behavior of  $\log \ell_{\alpha}$ along $\gamma$.  The second derivative is $(\ddot\ell_{\alpha}\ell_{\alpha}-\dot\ell_{\alpha}\dot\ell_{\alpha})\ell_{\alpha}^{-2}$.  Theorem \ref{conn} provides the expansion $2\ell_{\alpha}\ddot\ell_{\alpha}=(\dot\ell_{\alpha})^2+3(\dot\ell_{\alpha}\circ J)^2+O(\ell_{\alpha}^3\|U\|^2)$ and consequently the second derivative of $\log \ell_{\alpha}$ is given as 
$-2\ell_{\alpha}^{-1}\langle\lambda_{\alpha},U\rangle^2+6\ell_{\alpha}^{-1}\langle J\lambda_{\alpha},U\rangle^2+O(\ell_{\alpha}\|U\|^2)$.  From Theorem \ref{approx} the second derivative of $\log \ell_{\alpha}$ along $\gamma$ is bounded as $O(\ell_{\alpha}^{1/2}\|U\|^2)$.  In particular for $\gamma$ close to $\gamma_0$ the collar width is a shallow function along $\gamma$ with second derivative close to zero.  

We consider a lower bound for the weight function $\mathbb P_{\omega}$ of Section \ref{basics}.  Consider a marked surface $R$ with reflection $\upsilon$ and close to $R_0\cup_{p,q}\overline{R_0}$.  Introduce the reduced surface $R_{red}$ defined as the complement in $R$ of the unit area horoball neighborhoods of the cusps. The reduced surface $R_{red}$ is the union of thick regions, each containing a segment of $\omega$, each region with geometry close to a thick region of $R_0\cup_{p,q}\overline{R_0}$ and the collars $\mathcal C(p),\mathcal C(q)$, each crossed by a segment of $\omega$.  By elementary considerations the diameter of a thick region is bounded above in terms of a positive lower bound for the geodesic-lengths for $R_0$. The weight function is a sum of inverse square exponential distances to $\omega$ and so is bounded below on a thick region in terms of the diameter and is bounded below on a collar in terms of the maximal meridian length $2$.  In summary $\mathbb P_{\omega}$ is bounded below on $R_{red}$ in terms of a positive lower bound for the geodesic-lengths for $R_0$.  Following Theorem \ref{hessbd} for $\mu$ a harmonic Beltrami differential, we have $\int_{R_{red}}\mu\overline{\mu}\,dA\le c \hess\ell_{\omega}[\mu,\mu]$ for a constant depending on a positive lower bound for the geodesic-lengths for $R_0$.  At the end of Section 3.3 of \cite{Wlbhv} it is shown that the integral of $\mu\overline{\mu}$ over an area $1$ horoball is bounded above in terms of a universal constant and the integral over the complement of the area $1$ horoball in the area $2$ horoball.  In particular the integral $\int_{R_{red}}\mu\overline{\mu}\,dA$ is bounded below in terms of a universal positive constant and the WP pairing.    
The pairing for the tangent field of $\gamma$ converges to the pairing for $\gamma_0$.  In summary the convexity of $\ell_{\omega}$ along $\gamma$ is bounded below by a positive multiple of the WP pairing.  

In conclusion the function $\tilde{\ell}_{\omega_0}$ along $\gamma_0$ is a uniform limit of uniformly convex functions and so is strictly convex.  Finally for a cusp region the distance between given length horocycles is universal and so a change of choice of horocycle affects $\tilde{\ell}_{\omega_0}$ by adding a constant.  The proof is complete.


\end{document}